\documentclass[12pt]{article}
\usepackage{amssymb}
\newtheorem{theorem}{Theorem}[section]
\newtheorem{proposition}[theorem]{Proposition}
\newtheorem{lemma}[theorem]{Lemma}
\newtheorem{corollary}[theorem]{Corollary}
\def\bull{\vrule height .9ex width .8ex depth -.1ex}
\newenvironment{proof}{\smallbreak \noindent {\bf Proof.~}}
              {\unskip\nobreak\hfill\hskip 2em \bull\par\medbreak}
\newenvironment{proofof}[1]{\medbreak\noindent{\bf Proof of~#1.~}}
              {\unskip\nobreak\hfill\hskip 2em \bull\par\medbreak}
\def\bZ{\mathbb{Z}}
\def\cG{\mathcal{G}}

\def\si{\sigma}
\def\Aut{\mathop{\mathrm{Aut}}}
\def\id{\mathrm{id}}
\def\rev#1{\stackrel{\leftarrow}{#1}}
\newcounter{thesame}
\title{On a free group of transformations\\ defined by an automaton}
\author{Mariya Vorobets\thanks{%
        Partially supported by the NSF grants DMS-0308985 and DMS-0456185.}
        {} and Yaroslav Vorobets$^\fnsymbol{thesame}$\thanks{%
The second author is supported by a Clay Research Scholarship.}
}
\date{}
\begin{document}

\maketitle

\begin{abstract}
We prove that three automorphisms of the rooted binary tree defined by a
certain $3$-state automaton generate a free non-Abelian group of rank $3$.
\end{abstract}

\section{Introduction}\label{main}

A Mealy automaton $A$ over a finite alphabet $X$ is determined by the set
of internal states $Q$, the state transition function $\phi:Q\times X\to
Q$, and the output function $\psi:Q\times X\to X$.  The automaton starts
its work at some state $q\in Q$ and a sequence of letters $x_1,x_2,\dots,
x_n\in X$ is input into $A$.  The automaton uses the function $\psi$ to
produce an output sequence $y_1,y_2,\dots,y_n$ ($y_i\in X$) while changing
its states according to the function $\phi$.  This gives rise to a
transformation $A_q:X^*\to X^*$ of the set $X^*$ of finite words over
alphabet $X$.  A detailed account of the theory of Mealy automata is given
in the survey paper \cite{GNS}.

An automaton is called finite if it has only finitely many states.  An
automaton is called initial if it has a fixed initial state.  An initial
automaton over an alphabet $X$ defines a transformation of the set $X^*$.
A finite non-initial automaton defines a finite number of transformations,
each being assigned to an internal state.  Assuming all of them are
invertible, the finite automaton defines a finitely generated
transformation group.  This group is an example of a self-similar group
(see \cite{N}).

All invertible transformations defined by finite initial automata over a
fixed alphabet $X$ form a countable group $\cG(X)$ (see \cite{GNS}).  The
set $X^*$ is endowed with the structure of a rooted regular tree so that
elements of $\cG(X)$ are automorphisms of the tree.  The groups of finite
automaton transformations contain many finitely generated subgroups with
remarkable properties (e.g., infinite torsion groups, groups of
intermediate growth) but their complete structure is yet to be understood.

\begin{figure}[t]
%
\includegraphics{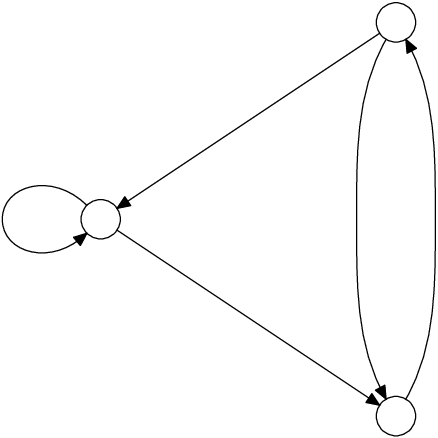}
\includegraphics{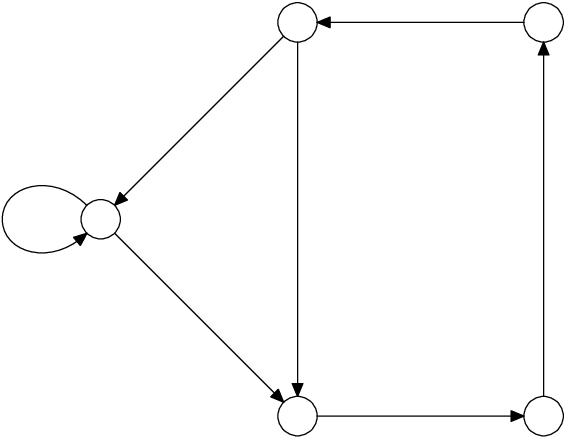}
\begin{picture}(0,0)(0,85)
\put(20,0)
{
 \begin{picture}(0,0)(0,0)
 \put(119,66){$a$}
 \put(34.5,-17){$b$}
 \put(119,-71){$c$}
 \put(66,33){\footnotesize $1|0$}
 \put(66,-39){\footnotesize $1|0$}
 \put(14,15){\footnotesize $0|1$}
 \put(95,-3){\footnotesize $0|1$}
 \put(131.5,-3){\footnotesize $\begin{array}{c} 0|0 \\ 1|1 \end{array}$}
 \end{picture}
}
\put(195.5,0)
{
 \begin{picture}(0,0)(0,0)
 \put(90.5,66){$a$}
 \put(34.5,-17){$b$}
 \put(90.5,-71){$c$}
 \put(161,-73){$d$}
 \put(162,66){$e$}
 \put(53,33){\footnotesize $1|0$}
 \put(53,-39){\footnotesize $1|0$}
 \put(14,15){\footnotesize $0|1$}
 \put(77.5,-3){\footnotesize $0|1$}
 \put(143,-3){\footnotesize $\begin{array}{c} 0|0 \\ 1|1 \end{array}$}
 \put(118,-45.5){\footnotesize $\begin{array}{c} 0|0 \\ 1|1 \end{array}$}
 \put(118,39.5){\footnotesize $\begin{array}{c} 0|0 \\ 1|1 \end{array}$}
 \end{picture}
}
\end{picture}
\vspace{165bp}
\caption{
\label{fig1}
Aleshin's automata.
}
\end{figure}

The problem of embedding the free non-Abelian group into a group $\cG(X)$
of finite automaton transformations turned out to be hard.  The first
attempt to solve it was made by Aleshin \cite{A}.  He introduced two finite
initial automata over alphabet $\{0,1\}$ and claimed that two
transformations of the rooted binary tree $\{0,1\}^*$ defined by these
automata generate a free group.  However the paper \cite{A} does not
contain a proof.  There is only a sketch that raises reasonable doubts as
to whether Aleshin's argument was complete.

Aleshin's automata are depicted in Figure \ref{fig1} by means of Moore
diagrams.  The Moore diagram of an automaton is a directed graph with
labeled edges.  The vertices are the states of the automaton and edges are
state transition routes.  Each label consists of two letters from the
alphabet.  The left one is the input field, it is used to choose a
transition route.  The right one is the output generated by the automaton.
Aleshin considered these automata as initial, with initial state $b$.

The first reliable result was obtained years later by Brunner and Sidki
\cite{BS}.  They showed that the group $\mathrm{GL}(n,\bZ)$ can be embedded
into the group of finite automaton transformations over the alphabet of
cardinality $2^n$.  Since $\mathrm{SL}(2,\bZ)$ contains a free group, this
yields two finite initial automata over a $4$-letter alphabet generating
the free group.  An example of a free group generated by finite automaton
transformations over a $2$-letter alphabet was given by Olijnyk and
Sushchanskij \cite{OS}.

A harder problem is to present the free group as the group defined by a
single finite non-initial automaton.  Recently it was solved by Glasner and
Mozes \cite{GM}.  They constructed infinitely many finite automata of
algebraic origin that define transformation groups with various properties,
in particular, free groups.  The simplest example from \cite{GM} is a pair
of automata $A_1$ and $A_2$.  The automaton $A_1$ is a $14$-state automaton
over a $6$-letter alphabet while $A_2$ is a $6$-state automaton over a
$14$-letter alphabet.  The automata define free groups on $7$ and $3$
generators, respectively.  Note that the number of free generators is half
of the number of states.  The $14$ transformations defined by the automaton
$A_1$ form a symmetric set of free generators.  That is, they split into
two sets of free generators such that the elements in each set are inverses
of elements in the other set.  The $6$ transformations defined by $A_2$
also form a symmetric set of free generators.  The automata $A_1$ and $A_2$
are dual in that either of them can be obtained from the other by
interchanging the alphabet with the set of internal states and the state
transition function with the output function.

In this paper, we consider the problem of finding a finite automaton that
defines a free group such that the rank of the group (i.e., the number of
free generators) is equal to the number of states of the automaton.  A
candidate for a solution has been available for some time.  Namely, Brunner
and Sidki conjectured (see \cite{S}) that the first of two Aleshin's
automata shown in Figure \ref{fig1} is the required one.  Here we prove
that this is indeed the case.

\begin{theorem}\label{main1}
Three automorphisms of the rooted binary tree defined
by the first Aleshin automaton generate a free transformation group on
three generators.
\end{theorem}

The proof of the theorem is based on the dual automaton approach (already
used in \cite{GM}).  The main idea of this approach is that properties of
the group defined by a finite automaton are determined by orbits of the
transformation group defined by the dual automaton.

The paper is organized as follows.  Section \ref{auto} addresses some
general constructions concerning automata and their properties.  In Section
\ref{free} we consider the Aleshin automaton along with a number of related
automata and establish some properties of these automata.  In Section
\ref{ind} we use results of Sections \ref{auto} and \ref{free} to prove
Theorem \ref{main1}.

\section{Automata}\label{auto}

An {\em automaton\/} $A$ is a quadruple $(Q,X,\phi,\psi)$ consisting of two
nonempty sets $Q$, $X$ along with two maps $\phi:Q\times X\to Q$,
$\psi:Q\times X\to X$.  The set $X$ is to be finite, it is called the {\em
input/output alphabet\/} of the automaton.  We say that $A$ is an automaton
over the alphabet $X$.  $Q$ is called the set of {\em internal states\/} of
$A$.  The automaton $A$ is called {\em finite\/} if the set $Q$ is finite.
$\phi$ and $\psi$ are called the {\em state transition function\/} and the
{\em output function}, respectively.  One may regard these functions as a
single map $(\phi,\psi):Q\times X\to Q\times X$.

The automaton $A$ functions as follows.  The active automaton is always
supposed to be in some state $q\in Q$.  The automaton reads an input letter
$x\in X$.  Then it performs two independent tasks: passes to the state
$\phi(q,x)$ and sends the letter $\psi(q,x)$ to the output.  After that the
automaton is ready to accept another input letter.  Usually the automaton
job consists of transducing the whole sequence of input letters.  Suppose
that $A$ started its work in a state $q\in Q$ (the {\em initial state\/})
and a word $w=x_1x_2\dots x_n$ over the alphabet $X$ was input into $A$.
As the result of automaton's job, we obtain two sequences: a sequence of
states $q_0=q,q_1,\dots,q_n$, which describes the internal work of the
automaton, and the output word $v=y_1y_2\dots y_n$.  Here
$q_i=\phi(q_{i-1},x_i)$ and $y_i=\psi(q_{i-1},x_i)$ for $1\le i\le n$.

Let $X^*$ denote the set of words over the alphabet $X$.  A {\em word\/}
$w\in X^*$ is merely a finite sequence whose elements belong to $X$.
However the elements of $w$ are called {\em letters\/} and $w$ is usually
written so that its elements are not separated by delimiters.  The number
of letters of a word $w$ is called its {\em length\/} and denoted by $|w|$.
It is assumed that $X^*$ contains the empty word $\varnothing$.  The set
$X$ is naturally embedded in $X^*$.  If $w_1=x_1\dots x_n$ and
$w_2=y_1\dots y_m$ are words over the alphabet $X$ then $w_1w_2$ denotes
their concatenation $x_1\dots x_ny_1\dots y_m$.  $X^*$ is the free monoid
generated by all elements of $X$ relative to the operation
$(w_1,w_2)\mapsto w_1w_2$.  Another structure on $X^*$ is that of a rooted
$k$-regular tree, where $k$ is the cardinality of $X$.  Namely, we consider
a graph with the set of vertices $X^*$ where two vertices $w_1,w_2\in X^*$
are joined by an edge if $w_1=w_2x$ or $w_2=w_1x$ for some $x\in X$.  The
root of the tree is the empty word.

In the above description of the automaton's functions, it is shown how the
output function $\psi$ can be extended to a map $\psi^*:Q\times X^*\to
X^*$.  By definition, $\psi^*$ satisfies the recursive relations
$\psi^*(q,\varnothing)=\varnothing$, $\psi^*(q,xw)=\psi(q,x)
\psi^*(\phi(q,x),w)$ for all $x\in X$, $w\in X^*$, $q\in Q$.  Moreover,
$\psi^*$ is uniquely determined by these relations.  Now for any $q\in Q$
we define a transformation $A_q:X^*\to X^*$, $A_q(w)=\psi^*(q,w)$ for all
$w\in X^*$.  We say that the transformation $A_q$ is defined by the
automaton $A$ with the initial state $q$.  Clearly, $A_q$ preserves length
of words.  Besides, $A_q$ transforms words from the left to the right, that
is, the first $n$ letters of a word $A_q(w)$ depend only on the first $n$
letters of $w$.  This implies that $A_q$ is an endomorphism of $X^*$ as a
rooted regular tree.  If $A_q$ is invertible then it belongs to the group
$\Aut(X^*)$ of automorphisms of the tree.

The semigroup of transformations of $X^*$ generated by $A_q$, $q\in Q$ is
denoted by $S(A)$.  The automaton $A$ is called {\em invertible\/} if $A_q$
is invertible for all $q\in Q$.  If $A$ is invertible then $A_q$, $q\in Q$
generate a transformation group $G(A)$, which is a subgroup of $\Aut(X^*)$.
We say that $S(A)$ (resp. $G(A)$) is the semigroup (resp. group) defined by
the automaton $A$.

\begin{lemma}\label{auto1}
Suppose the automaton $A$ is invertible.  Then the actions of the semigroup
$S(A)$ and the group $G(A)$ on $X^*$ have the same orbits.
\end{lemma}

\begin{proof}
The action of the group $G(A)$ on $X^*$ preserves length of words.  Hence
an arbitrary word $w\in X^*$ belongs to some finite set $W\subset X^*$ that
is invariant under the $G(A)$ action.  Since the automaton $A$ is
invertible, every transformation $g\in S(A)$ is invertible.  The
restriction of $g$ to $W$ is a bijection, which is of finite order since
the set $W$ is finite.  It follows that the semigroup $S(A)$ of
transformations of $X^*$ becomes a group when restricted to $W$.  In
particular, the restrictions to $W$ of the semigroup $S(A)$ and of the
group $G(A)$ coincide.  Then the orbits $\{g(w)\mid g\in S(A)\}$ and
$\{g(w)\mid g\in G(A)\}$ coincide as well.
\end{proof}

One way to picture an automaton, which we use in this paper, is the {\em
Moore diagram}.  The Moore diagram of an automaton $A=(Q,X,\phi,\psi)$ is
a directed graph with labeled edges defined as follows.  The vertices of
the graph are states of the automaton $A$.  Every edge carries a
label of the form $x|y$, where $x,y\in X$.  The left field $x$ of the label
is referred to as the {\em input field\/} while the right field $y$ is
referred to as the {\em output field}.  The set of edges of the graph is in
a one-to-one correspondence with the set $Q\times X$.  Namely, for any
$q\in Q$ and $x\in X$ there is an edge that goes from the vertex $q$ to
$\phi(q,x)$ and carries the label $x|\psi(q,x)$.  The Moore diagram of an
automaton can have loops (edges joining a vertex to itself) and multiple
edges.  To simplify pictures, we do not draw multiple edges in this paper.
Instead, we use multiple labels.

The transformations $A_q$, $q\in Q$ can be defined in terms of the Moore
diagram of the automaton $A$.  For any $q\in Q$ and $w\in X^*$ we find a
path $\delta$ in the Moore diagram such that $\delta$ starts at the vertex
$q$ and the word $w$ can be obtained by reading the input fields of labels
along $\delta$.  Such a path exists and is unique.  Then the word $A_q(w)$
is obtained by reading the output fields of labels along the path $\delta$.

Let $\Gamma$ denote the Moore diagram of the automaton $A$.  We associate
to $\Gamma$ two directed graphs $\Gamma_1$ and $\Gamma_2$ with labeled
edges.  $\Gamma_1$ is obtained from $\Gamma$ by interchanging the input and
output fields of all labels.  That is, a label $x|y$ is replaced by $y|x$.
$\Gamma_2$ is obtained from $\Gamma$ by reversing all edges.
The {\em inverse automaton\/} of $A$ is the automaton whose Moore diagram
is $\Gamma_1$.  The {\em reverse automaton\/} of $A$ is the automaton whose
Moore diagram is $\Gamma_2$.  If one of the graphs $\Gamma_1$ and
$\Gamma_2$ is the Moore diagram of an automaton then the automaton shares
with $A$ the alphabet and internal states.  Moreover, its state transition
and output functions are uniquely determined by the graph.  However neither
$\Gamma_1$ nor $\Gamma_2$ must be the Moore diagram of an automaton.  The
necessary and sufficient condition for such a labelling to be the Moore
diagram of an automaton is that for any $q\in Q$ and $x\in X$ there is
exactly one edge of the graph that goes out of the vertex $q$ and has $x$
as the input field of its label.

\begin{proposition}\label{auto2}
An automaton $A=(Q,X,\phi,\psi)$ is invertible if and only if for any $q\in
Q$ the map $\psi(q,\cdot):X\to X$ is bijective.  The inverse automaton $I$
is well defined if and only if $A$ is invertible.  If this is the case,
then $I_q=A_q^{-1}$ for all $q\in Q$.
\end{proposition}

\begin{proof}
Let $\Gamma$ be the Moore diagram of $A$ and $\Gamma_1$ be the graph
obtained from $\Gamma$ by interchanging the input and output fields of all
labels.  For any $x\in X$ and $q\in Q$ the number of edges of $\Gamma_1$
that go out of the vertex $q$ and have $x$ as the input field of their
labels is equal to the number of $y\in X$ such that $\psi(y,q)=x$.  The
inverse automaton of $A$ is well defined if this number is always equal to
$1$, that is, if and only if the map $\psi(\cdot,q):X\to X$ is bijective
for any $q\in Q$.

Suppose there exists $q\in Q$ such that the map $\psi(\cdot,q)$ is not
bijective.  Since $X$ is a finite set, $\psi(x_1,q)=\psi(x_2,q)$ for some
$x_1,x_2\in X$, $x_1\ne x_2$.  Then $A_q(x_1)=A_q(x_2)$; in particular,
$A_q$ is not invertible.

Now assume that the map $\psi(\cdot,q)$ is bijective for any $q\in Q$.  By
the above the inverse automaton $I$ of $A$ is well defined.  The graphs
$\Gamma$ and $\Gamma_1$, which are the Moore diagrams of automata $A$ and
$I$, differ only in labels.  In particular, they share all paths.  Given
such a path $\delta$, a word $w\in X^*$ can be obtained by reading the
input fields of labels along $\delta$ in one of the graphs if and only if
the same word $w$ is obtained by reading the output fields of labels along
$\delta$ in the other graph.  It follows that for any $q\in Q$ and
$w_1,w_2\in X^*$ we have $w_2=A_q(w_1)$ if and only if $w_1=I_q(w_2)$.
This means that $A_q$ is invertible for any $q\in Q$ and $I_q=A_q^{-1}$.
\end{proof}

For any word $w=x_1x_2\dots x_n$ over an alphabet $X$ we denote by
$\rev{w}$ the reversed word $x_n\dots x_2x_1$.

\begin{proposition}\label{auto3}
Given an automaton $A=(Q,X,\phi,\psi)$, the reverse automaton $R$ is well
defined if and only if for any $x\in X$ the map $\phi(\cdot,x):Q\to Q$ is
bijective.  Assume this is the case and let $w_1,w_2\in X^*$.  Then
$w_2=g(w_1)$ for some $g\in S(A)$ if and only if $\rev{w_2}=h(\rev{w_1})$
for some $h\in S(R)$.
\end{proposition}

\begin{proof}
Let $\Gamma$ be the Moore diagram of $A$ and $\Gamma_2$ be the graph
obtained from $\Gamma$ by reversing all edges.  For any $q\in Q$ and $x\in
X$ the number of edges of $\Gamma_2$ that go out of the vertex $q$ and have
$x$ as the input field of their labels is equal to the number of $p\in Q$
such that $\phi(p,x)=q$.  The reverse automaton $R$ of $A$ is well
defined if this number is always equal to $1$, that is, if and only if
$\phi(\cdot,x):Q\to Q$ is bijective for any $x\in X$.

Assume the reverse automaton $R$ is well defined.  The graph $\Gamma_2$ is
the Moore diagram of $R$.  Any path $\delta$ in the graph $\Gamma$ is
assigned a ``reversed'' path $\delta'$ in the graph $\Gamma_2$ such that a
word $w\in X^*$ can be obtained by reading the input (resp. output) fields
of labels along $\delta$ if and only if the reversed word $\rev{w}$ is
obtained by reading the input (resp. output) fields of labels along
$\delta'$.  It follows that, given $w_1,w_2\in X^*$, we have $w_2=A_q(w_1)$
for some $q\in Q$ if and only if $\rev{w_2}=R_p(\rev{w_1})$ for some $p\in
Q$.

Let $w_1,w_2\in X^*$.  If $w_2=g(w_1)$ for some $g\in S(A)$ then there is a
sequence of words $u_1=w_1,u_2,\dots,u_n=w_2$ ($n>1$) such that
$u_{i+1}=A_{q_i}(u_i)$ for some $q_i\in Q$, $1\le i\le n-1$.  By the above
$\rev{u_{i+1}}=R_{p_i}(\rev{u_i})$, where $p_i\in Q$, $1\le i\le n-1$.
Hence $\rev{w_2}=h(\rev{w_1})$ for some $h\in S(R)$.  Similarly,
$\rev{w_2}=h(\rev{w_1})$ for some $h\in S(R)$ implies that $w_2=g(w_1)$ for
some $g\in S(A)$.
\end{proof}

Let $A=(Q,X,\phi,\psi)$ be an automaton.  For any nonempty word
$\xi=q_1q_2\dots q_n\in Q^*$ we let $A_\xi=A_{q_n}\dots A_{q_2}A_{q_1}$.
Also, we let $A_\varnothing=1$ (here $1$ stands for the unit element of the
group $\Aut(X^*)$, i.e., the identity mapping on $X^*$).  Obviously, any
element of the semigroup $S(A)$ is represented as $A_\xi$ for a nonempty
word $\xi\in Q^*$.  The map $X^*\times Q^*\to X^*$ given by $(w,\xi)\mapsto
A_\xi(w)$ defines a right action of the monoid $Q^*$ on the rooted regular
tree $X^*$.  That is, $A_{\xi_1\xi_2}(w)=A_{\xi_2}(A_{\xi_1}(w))$ for all
$\xi_1,\xi_2\in Q^*$ and $w\in X^*$.

To each finite automaton $A=(Q,X,\phi,\psi)$ we associate a {\em dual
automaton\/} $D$, which is obtained from $A$ by interchanging the alphabet
with the set of internal states and the state transition function with the
output function.  To be precise, $D=(X,Q,\tilde\phi,\tilde\psi)$, where
$\tilde\phi(x,q)=\psi(q,x)$ and $\tilde\psi(x,q)=\phi(q,x)$ for all $x\in
X$ and $q\in Q$.  Unlike the inverse and reverse automata, the dual
automaton is always well defined.  It is easy to see that $A$ is the dual
automaton of $D$.

The dual automaton $D$ defines a right action of the monoid $X^*$ on $Q^*$
given by $(\xi,w)\mapsto D_w(\xi)$.  This action and the action of $Q^*$ on
$X^*$ defined by the automaton $A$ are related in the following way.

\begin{proposition}\label{auto4}
For any $w,u\in X^*$ and $\xi,\eta\in Q^*$,
$$
A_\xi(wu)=A_\xi(w)A_{D_w(\xi)}(u),
$$ $$
D_w(\xi\eta)=D_w(\xi)D_{A_\xi(w)}(\eta).
$$
\end{proposition}

\begin{proof}
First consider the case when one of words $w$ and $\xi$ is empty.  If
$w=\varnothing$ then $A_\xi(w)=\varnothing$ and $D_w=1$.  If
$\xi=\varnothing$ then $A_\xi=1$ and $D_w(\xi)=\varnothing$.

Further, consider the case when both $w$ and $\xi$ are one-letter words.
In this case, $A_\xi(wu)=xA_q(u)$, where $x=\psi(\xi,w)$, $q=\phi(\xi,w)$,
while $D_w(\xi\eta)=pD_y(\eta)$, where $p=\tilde\psi(w,\xi)$,
$y=\tilde\phi(w,\xi)$.  It remains to observe that $\tilde\phi(w,\xi)=
\psi(\xi,w)=A_\xi(w)$ and $\phi(\xi,w)=\tilde\psi(w,\xi)=D_w(\xi)$.

In the general case, we prove the proposition by induction on the sum of
lengths of words $w$ and $\xi$.  By the above the proposition holds when
$|w|+|\xi|\le2$.  Now let $n>2$ and assume the proposition holds if
$|w|+|\xi|<n$.  Take arbitrary words $w,u\in X^*$ and $\xi,\eta\in Q^*$
such that $|w|+|\xi|=n$.  Since $n>2$, at least one of the words $w$, $\xi$
has more than one letter.   If $|w|>1$ then $w=w_1w_2$ for some words
$w_1,w_2\in X^*$ that are shorter than $w$.  Repeatedly using the inductive
assumption, we obtain
$$
A_\xi(wu)=A_\xi(w_1w_2u)=A_\xi(w_1)A_{D_{w_1}(\xi)}(w_2u)=
$$ $$
A_\xi(w_1)A_{D_{w_1}(\xi)}(w_2)A_{D_{w_2}(D_{w_1}(\xi))}(u)=
A_\xi(w_1w_2)A_{D_{w_1w_2}(\xi)}(u).
$$
If $|\xi|>1$ then $\xi=\xi_1\xi_2$ for some words $\xi_1,\xi_2\in Q^*$ that
are shorter than $\xi$.  By the inductive assumption,
$$
A_\xi(wu)=A_{\xi_1\xi_2}(wu)=A_{\xi_2}(A_{\xi_1}(wu))=
A_{\xi_2}(A_{\xi_1}(w)A_{D_w(\xi_1)}(u))=
$$ $$
A_{\xi_2}(A_{\xi_1}(w))A_\zeta(A_{D_w(\xi_1)}(u))=
A_{\xi_1\xi_2}(w)A_{D_w(\xi_1)\zeta}(u),
$$
where $\zeta=D_{A_{\xi_1}(w)}(\xi_2)$.  Since $|w|+|\xi_1|<n$, we have
$D_w(\xi_1)\zeta=D_w(\xi_1\xi_2)$.

The equality $D_w(\xi\eta)=D_w(\xi)D_{A_\xi(w)}(\eta)$ is verified in the
same way.
\end{proof}

\begin{corollary}\label{auto5}
Suppose $A_\xi=1$ for some $\xi\in Q^*$.  Then $A_{g(\xi)}=1$ for every
$g\in S(D)$.
\end{corollary}

\begin{proof}
$g\in S(D)$ means that $g=D_w$ for some word $w\in X^*$.  By Proposition
\ref{auto4}, $A_\xi(wu)=A_\xi(w)A_{g(\xi)}(u)$ for any $u\in X^*$.  Since
$A_\xi=1$, we have $wu=wA_{g(\xi)}(u)$, which implies that
$A_{g(\xi)}(u)=u$.
\end{proof}

Let $A=(Q_1,X,\phi_1,\psi_1)$ and $B=(Q_2,X,\phi_2,\psi_2)$ be two automata
over the same alphabet $X$.  Assume $Q_1\cap Q_2=\emptyset$.  The {\em
disjoint union\/} of automata $A$ and $B$ is the automaton $U=(Q_1\cup Q_2,
X,\phi,\psi)$, where functions $\phi$, $\psi$ are such that $\phi=\phi_1$
and $\psi=\psi_1$ on $Q_1\times X$ while $\phi=\phi_2$ and $\psi=\psi_2$ on
$Q_2\times X$.  Obviously, $U_q=A_q$ for $q\in Q_1$ and $U_q=B_q$ for $q\in
Q_2$.  The Moore diagram of $U$ is the disjoint union of the Moore diagrams
of $A$ and $B$.  For example, Figure \ref{fig2} may be regarded either as
the Moore diagrams of two $3$-state automata or as the Moore diagram of a
single $6$-state automaton, their disjoint union.  If $Q_1\cap Q_2\ne
\emptyset$, then the disjoint union of automata $A$ and $B$ is not defined
(this is the case for two automata whose Moore diagrams are depicted in
Figure \ref{fig1}).  However we can rename some states of the automata so
that they share no states anymore and then consider the disjoint union of
the modified automata.

\section{The Aleshin automaton}\label{free}

In this section we turn to the main object of study of this paper, which is
the first of the two automata of Aleshin whose Moore diagrams are depicted
in Figure \ref{fig1}.  This is an automaton $A$ over the alphabet
$X=\{0,1\}$ with the set of internal states $Q=\{a,b,c\}$.  The state
transition function $\phi$ and the output function $\psi$ of $A$ are
defined as follows: $\phi(a,0)=\phi(b,1)=c$, $\phi(a,1)=\phi(b,0)=b$,
$\phi(c,0)=\phi(c,1)=a$; $\psi(a,0)=\psi(b,0)=\psi(c,1)=1$, $\psi(a,1)=
\psi(b,1)=\psi(c,0)=0$.

\begin{figure}[t]
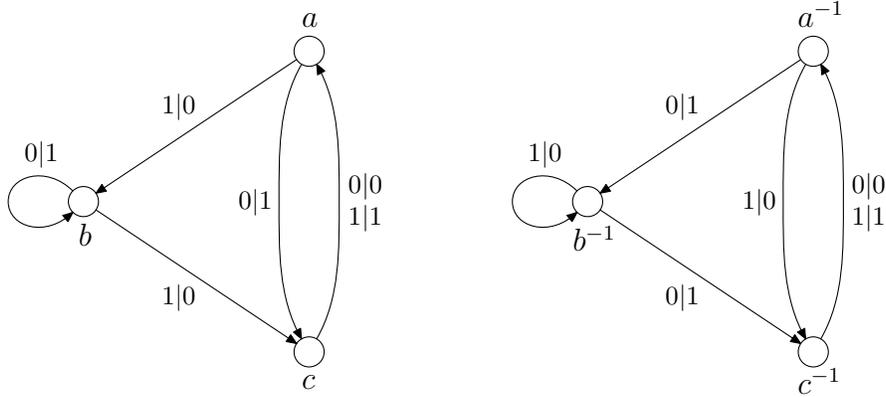

%
\includegraphics{first.eps}
\includegraphics{first.eps}
\begin{picture}(0,0)(0,80)
\put(20,0)
{
 \begin{picture}(0,0)(0,0)
 \put(119,66){$a$}
 \put(34.5,-17){$b$}
 \put(119,-71){$c$}
 \put(66,33){\footnotesize $1|0$}
 \put(66,-39){\footnotesize $1|0$}
 \put(14,15){\footnotesize $0|1$}
 \put(95,-3){\footnotesize $0|1$}
 \put(131.5,-3){\footnotesize $\begin{array}{c} 0|0 \\ 1|1 \end{array}$}
 \end{picture}
}
\put(210.5,0)
{
 \begin{picture}(0,0)(0,0)
 \put(116,66){$a^{-1}$}
 \put(31,-19){$b^{-1}$}
 \put(116,-73){$c^{-1}$}
 \put(66,33){\footnotesize $0|1$}
 \put(66,-39){\footnotesize $0|1$}
 \put(14,15){\footnotesize $1|0$}
 \put(95,-3){\footnotesize $1|0$}
 \put(131.5,-3){\footnotesize $\begin{array}{c} 0|0 \\ 1|1 \end{array}$}
 \end{picture}
}
\end{picture}
\vspace{155bp}
\caption{
\label{fig2}
The Aleshin automaton and its inverse.
}
\end{figure}

Using Proposition \ref{auto2}, it is easy to verify that $A$ is an
invertible automaton.  In particular, $A_a$, $A_b$, and $A_c$ are
automorphisms of the rooted binary tree $X^*$.  Let $I$ denote the inverse
automaton of $A$.  By $I'$ denote the automaton obtained from $I$ by
renaming its states $a$, $b$, $c$ to $a^{-1}$, $b^{-1}$, $c^{-1}$,
respectively.  Here, $a^{-1}$, $b^{-1}$, and $c^{-1}$ are regarded as
elements of the free group on generators $a$, $b$, $c$.  Further, let $B$
denote the disjoint union of automata $A$ and $I'$.  $B$ is an automaton
over the alphabet $X=\{0,1\}$ with the set of internal states
$Q_\pm=\{a,b,c,a^{-1},b^{-1},c^{-1}\}$.  We have $B_a=A_a$, $B_b=A_b$,
$B_c=A_c$, $B_{a^{-1}}=A_a^{-1}$, $B_{b^{-1}}=A_b^{-1}$, $B_{c^{-1}}=
A_c^{-1}$.  In particular, $S(B)=G(A)$.  Figure \ref{fig2} depicts the
Moore diagram of the automaton $B$.

\begin{figure}[b]
%
\includegraphics{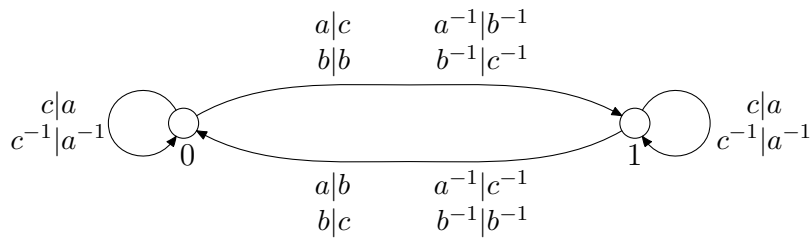}
\begin{picture}(0,0)(-100,55)
\put(-1,-16){$0$}
\put(168,-16){$1$}
\put(45,27){\small $\begin{array}{r} a|c \\ b|b \end{array}$}
\put(90,27){\small $\begin{array}{r} a^{-1}|b^{-1} \\ b^{-1}|c^{-1} \end{array}$}
\put(45,-34){\small $\begin{array}{r} a|b \\ b|c \end{array}$}
\put(90,-34){\small $\begin{array}{r} a^{-1}|c^{-1} \\ b^{-1}|b^{-1} \end{array}$}
\put(-70,-3){\small $\begin{array}{c} c|a \\ c^{-1}|a^{-1} \end{array}$}
\put(197,-3){\small $\begin{array}{c} c|a \\ c^{-1}|a^{-1} \end{array}$}
\end{picture}
\vspace{105bp}
\caption{
\label{fig3}
The dual automaton $D$.
}
\end{figure}

\begin{figure}[tb]
%
\includegraphics{dual.eps}
\includegraphics{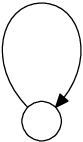}
\begin{picture}(0,0)(-50,55)
\put(-2,-14){$\alpha$}
\put(165,-16){$\beta$}
\put(45,27){\small $\begin{array}{c} a|a \\ b|b \end{array}$}
\put(90,27){\small $\begin{array}{c} a^{-1}|b^{-1} \\ \,b^{-1}|a^{-1} \end{array}$}
\put(45,-34){\small $\begin{array}{c} a|b \\ \,b|a \end{array}$}
\put(90,-34){\small $\begin{array}{c} a^{-1}|a^{-1} \\ b^{-1}|b^{-1} \end{array}$}
\put(-38,25){\small $\begin{array}{c} c|c \\ c^{-1}|c^{-1} \end{array}$}
\put(170,25){\small $\begin{array}{c} c|c \\ c^{-1}|c^{-1} \end{array}$}
\put(267,-14){$\gamma$}
\put(238,17){\small $\begin{array}{c} a|a \\ b|c \\ c|b \end{array}$}
\put(280,17){\small $\begin{array}{c} a^{-1}|a^{-1} \\ b^{-1}|c^{-1} \\ c^{-1}|b^{-1} \end{array}$}
\end{picture}
\vspace{105bp}
\caption{
\label{fig4}
Automaton $E$.
}
\end{figure}

We are going to consider two automata over the alphabet $Q_\pm$.  The first
one is the dual automaton $D$ of the automaton $B$.  It has two states $0$
and $1$; its Moore diagram is shown in Figure \ref{fig3}.  The other one is
an auxiliary automaton $E$.  The automaton $E$ has $3$ internal states
$\alpha$, $\beta$, and $\gamma$.  Its transition function $\phi_E$ is
defined as follows.  If $q\in\{a,b,a^{-1},b^{-1}\}$ then $\phi_E(\alpha,q)=
\beta$ and $\phi_E(\beta,q)=\alpha$.  If $q\in\{c,c^{-1}\}$ then
$\phi_E(\alpha,q)=\alpha$ and $\phi_E(\beta,q)=\beta$.  Also,
$\phi_E(\gamma,q)=\gamma$ for all $q\in Q_\pm$.  The output function
$\psi_E$ of $E$ is defined so that $\psi_E(\alpha,q)=\si_\alpha(q)$,
$\psi_E(\beta,q)=\si_\beta(q)$, and $\psi_E(\gamma,q)=\si_\gamma(q)$ for
all $q\in Q_\pm$, where $\si_\alpha=(a^{-1}b^{-1})$, $\si_\beta=(ab)$, and
$\si_\gamma=(bc)(b^{-1}c^{-1})$ are permutations on the set $Q_\pm$.  The
Moore diagram of the automaton $E$ is shown in Figure \ref{fig4}.  Using
Proposition \ref{auto2}, one can verify that $D$ and $E$ are invertible
automata.

To each permutation $\tau$ on the set $Q=\{a,b,c\}$ we assign a mapping
$\pi_\tau:Q_\pm^*\to Q_\pm^*$ as follows.  First we define a permutation
$\tilde\tau$ on the set $Q_\pm$ by $\tilde\tau(q)=\tau(q)$,
$\tilde\tau(q^{-1})=(\tau(q))^{-1}$ for all $q\in Q$.  Then for any
nonempty word $\xi=q_1q_2\dots q_n$ over the alphabet $Q_\pm$ we let
$\pi_\tau(\xi)=\tilde\tau(q_1)\tilde\tau(q_2)\dots\tilde\tau(q_n)$.
Besides, we let $\pi_\tau(\varnothing)=\varnothing$.  Clearly, the mapping
$\pi_\tau$ is an automorphism of the monoid $Q_\pm^*$.

\begin{lemma}\label{free1}
$\pi_{\tau_1\tau_2}=\pi_{\tau_1}\pi_{\tau_2}$ and $\pi_{\tau^{-1}}=
\pi_\tau^{-1}$ for any permutations $\tau,\tau_1,\tau_2$ on the set
$\{a,b,c\}$.
\end{lemma}

\begin{proof}
Let $\tau_1$, $\tau_2$ be permutations on $\{a,b,c\}$.  Since
$\pi_{\tau_1}$, $\pi_{\tau_2}$, and $\pi_{\tau_1\tau_2}$ are automorphisms
of the monoid $Q_\pm^*$, it follows that $\pi_{\tau_1\tau_2}=\pi_{\tau_1}
\pi_{\tau_2}$ whenever $\pi_{\tau_1\tau_2}(q)=\pi_{\tau_1}\pi_{\tau_2}(q)$
for all $q\in Q_\pm$.  Given $q\in\{a,b,c\}$, we have
$\pi_{\tau_1}\pi_{\tau_2}(q)=\pi_{\tau_1}(\tau_2(q))=\tau_1\tau_2(q)=
\pi_{\tau_1\tau_2}(q)$, $\pi_{\tau_1}\pi_{\tau_2}(q^{-1})=
\pi_{\tau_1}((\tau_2(q))^{-1})=(\tau_1\tau_2(q))^{-1}=
\pi_{\tau_1\tau_2}(q^{-1})$.

Let $\tau$ be a permutation on $\{a,b,c\}$.  By the above
$\pi_\tau\pi_{\tau^{-1}}=\pi_{\tau^{-1}}\pi_\tau=\pi_{\id}$, where $\id$
denotes the identity map on $\{a,b,c\}$.  Obviously, $\pi_{\id}$ is the
identity map on $Q_\pm^*$.  Hence $\pi_{\tau^{-1}}=\pi_\tau^{-1}$.
\end{proof}

\begin{lemma}\label{free2}
$E_\alpha^2=E_\beta^2=E_\gamma^2=1$, $E_\alpha E_\beta=E_\beta E_\alpha=
\pi_{(ab)}$, $E_\gamma=\pi_{(bc)}$.

The group $G(E)$ contains $\pi_\tau$ for any permutation $\tau$ on
$\{a,b,c\}$.
\end{lemma}

\begin{proof}
It is easy to see that the inverse automaton of $E$ coincides with $E$.
Proposition \ref{auto2} implies that $E_\alpha^2=E_\beta^2=E_\gamma^2=1$.
The equality $E_\gamma=\pi_{(bc)}$ follows from the definition of the
automaton $E$.

Consider the following permutations on the set $Q_\pm$: $\si_\alpha=
(a^{-1}b^{-1})$, $\si_\beta=(ab)$, and $\si=(ab)(a^{-1}b^{-1})$.  Given an
arbitrary word $\xi=q_1q_2\dots q_k\in Q_\pm^*$, the words $E_\alpha(\xi)=
a_1a_2\dots a_k$, $E_\beta(\xi)=b_1b_2\dots b_k$, and $\pi_{(ab)}(\xi)=
c_1c_2\dots c_k$ are computed as follows.  For any $i$, $1\le i\le k$ we
count the number of times when letters $a,b,a^{-1},b^{-1}$ occur in the
sequence $q_1,\dots,q_{i-1}$.  If this number is even then $a_i=
\si_\alpha(q_i)$, $b_i=\si_\beta(q_i)$; otherwise $a_i=\si_\beta(q_i)$,
$b_i=\si_\alpha(q_i)$.  In any case, $c_i=\si(q_i)$.  Since the set
$\{a,b,a^{-1},b^{-1}\}$ is invariant under permutations $\si_\alpha$ and
$\si_\beta$, the equalities $\si_\alpha\si_\beta=\si_\beta\si_\alpha=\si$
imply that $E_\alpha E_\beta=E_\beta E_\alpha=\pi_{(ab)}$.

Since $\pi_{(ab)},\pi_{(bc)}\in G(E)$ and the group of permutations on
$\{a,b,c\}$ is generated by permutations $(ab)$ and $(bc)$, it follows from
Lemma \ref{free1} that $G(E)$ contains all transformations of the form
$\pi_\tau$.
\end{proof}

\begin{lemma}\label{free3}
$D_0=\pi_{(ac)}E_\alpha=\pi_{(abc)}E_\beta$, $D_1=\pi_{(abc)}E_\alpha=
\pi_{(ac)}E_\beta$.
\end{lemma}

\begin{proof}
Let $\xi=q_1q_2\dots q_k$ be a word in the alphabet $Q_\pm$.  According to
the definitions of automata $D$ and $E$, the words $E_\alpha(\xi)=
a_1a_2\dots a_k$, $E_\beta(\xi)=b_1b_2\dots b_k$, $D_0(\xi)=c_1c_2\dots
c_k$, and $D_1(\xi)=d_1d_2\dots d_k$ can be computed as follows.  First we
define $4$ permutations on the set $Q_\pm$: $\si_\alpha=(a^{-1}b^{-1})$,
$\si_\beta=(ab)$, $\si_0=(ac)(a^{-1}b^{-1}c^{-1})$, $\si_1=
(abc)(a^{-1}c^{-1})$.  Now for any $i$, $1\le i\le k$ we count the number
of times when letters $a,b,a^{-1},b^{-1}$ occur in the sequence
$q_1,\dots,q_{i-1}$.  If this number is even then $a_i=\si_\alpha(q_i)$,
$b_i=\si_\beta(q_i)$, $c_i=\si_0(q_i)$, $d_i=\si_1(q_i)$.  Otherwise $a_i=
\si_\beta(q_i)$, $b_i=\si_\alpha(q_i)$, $c_i=\si_1(q_i)$, $d_i=\si_0(q_i)$.

Let us consider $2$ more permutations on the set $Q_\pm$: $\tau_0=
(ac)(a^{-1}c^{-1})$ and $\tau_1=(abc)(a^{-1}b^{-1}c^{-1})$.  It is easy to
verify that $\si_0=\tau_0\si_\alpha=\tau_1\si_\beta$ and $\si_1=
\tau_1\si_\alpha=\tau_0\si_\beta$.  Therefore for any $i$, $1\le i\le k$ we
have $c_i=\tau_0(a_i)=\tau_1(b_i)$ and $d_i=\tau_1(a_i)=\tau_0(b_i)$.  This
means that $D_0(\xi)=\pi_{(ac)}E_\alpha(\xi)=\pi_{(abc)}E_\beta(\xi)$ and
$D_1(\xi)=\pi_{(abc)}E_\alpha(\xi)=\pi_{(ac)}E_\beta(\xi)$.
\end{proof}

\begin{proposition}\label{free4}
$G(D)=G(E)$.
\end{proposition}

\begin{proof}
Lemmas \ref{free2} and \ref{free3} imply that $D_0,D_1\in G(E)$; therefore
$G(D)\subset G(E)$.

By Lemma \ref{free3}, $D_0D_1^{-1}=\pi_{(ac)}E_\alpha
(\pi_{(abc)}E_\alpha)^{-1}=\pi_{(ac)}\pi_{(abc)}^{-1}$.  By Lemma
\ref{free1}, $\pi_{(ac)}\pi_{(abc)}^{-1}=\pi_{(ac)}\pi_{(acb)}=\pi_{(bc)}=
E_\gamma$.  Similarly,
$$
D_0^{-1}D_1=(\pi_{(ac)}E_\alpha)^{-1}\pi_{(abc)}E_\alpha=E_\alpha^{-1}
\pi_{(ac)}^{-1}\pi_{(abc)}E_\alpha=E_\alpha^{-1}\pi_{(ab)}E_\alpha.
$$
Lemma \ref{free2} implies that $E_\alpha$ and $\pi_{(ab)}$ commute, hence
$D_0^{-1}D_1=\pi_{(ab)}$.  The group of all permutations on the set
$\{a,b,c\}$ is generated by permutations $(ab)$ and $(bc)$.  Since
$\pi_{(ab)},\pi_{(bc)}\in G(D)$, it follows from Lemma \ref{free1} that
$G(D)$ contains all transformations of the form $\pi_\tau$.  Then
$E_\alpha,E_\beta\in G(D)$ by Lemma \ref{free3}.  Thus $G(E)\subset G(D)$.
\end{proof}

\section{Patterns and orbits}\label{ind}

This section is devoted to the proof of Theorem \ref{main1}.  We use the
notation of the previous section.

As shown in Section \ref{auto}, the automaton $B$ defines a right action
$X^*\times Q_\pm^*\to X^*$ of the monoid $Q_\pm^*$ on the rooted binary
tree $X^*$ given by $(w,\xi)\mapsto B_\xi(w)$.  Since $Q_\pm^*$ is the free
monoid generated by $a,b,c,a^{-1},b^{-1},c^{-1}$, there exists a unique
homomorphism $\chi:Q_\pm^*\to\{-1,1\}$ such that $\chi(a)=\chi(b)=
\chi(a^{-1})=\chi(b^{-1})=-1$, $\chi(c)=\chi(c^{-1})=1$.

\begin{lemma}\label{ind1}
Given $\xi\in Q_\pm^*$, the automorphism $B_\xi$ of the rooted binary tree
$\{0,1\}^*$ acts trivially on the first level of the tree (i.e., on
one-letter words) if and only if $\chi(\xi)=1$.
\end{lemma}

\begin{proof}
For any $\xi\in Q_\pm^*$ let $\tilde\chi(\xi)=1$ if $B_\xi$ acts trivially
on the first level of the binary tree $X^*$ and $\tilde\chi(\xi)=-1$
otherwise.  If $\tilde\chi(\xi)=-1$ then $B_\xi(0)=1$, $B_\xi(1)=0$.  Since
$B_{\xi_1\xi_2}=B_{\xi_2}B_{\xi_1}$ for all $\xi_1,\xi_2\in Q_\pm^*$, it
follows that the map $\tilde\chi:Q_\pm^*\to\{-1,1\}$ is a homomorphism of
the monoid $Q_\pm^*$.  By definition, $\tilde\chi(a)=\tilde\chi(b)=-1$,
$\tilde\chi(c)=1$, and $\tilde\chi(q^{-1})=\tilde\chi(q)$ for any
$q\in\{a,b,c\}$.  Since $\tilde\chi$ is a homomorphism, it follows that
$\tilde\chi=\chi$.
\end{proof}

Now we introduce an alphabet consisting of two symbols $*$ and $*^{-1}$.  A
word over the alphabet $\{*,*^{-1}\}$ is called a {\em pattern}.  Every
word $\xi$ over the alphabet $Q_\pm$ is assigned a pattern $v$ that is
obtained from $\xi$ by substituting $*$ for each occurrence of letters
$a,b,c$ and substituting $*^{-1}$ for each occurrence of letters
$a^{-1},b^{-1},c^{-1}$.  We say that $v$ is the pattern of $\xi$ or that
$\xi$ follows the pattern $v$.

A word $\xi=q_1q_2\dots q_n\in Q_\pm^*$ is called {\em freely
irreducible\/} if none of its two-letter subwords $q_1q_2,q_2q_3,\dots,
q_{n-1}q_n$ coincides with one of the following words: $aa^{-1},bb^{-1},
cc^{-1},a^{-1}a,b^{-1}b,c^{-1}c$.  Otherwise $\xi$ is called {\em freely
reducible}.

\begin{lemma}\label{ind2}
For any nonempty pattern $v$ there exist words $\xi_1,\xi_2\in Q_\pm^*$
such that $\xi_1$ and $\xi_2$ are freely irreducible, follow the pattern
$v$, and $\chi(\xi_2)=-\chi(\xi_1)$.
\end{lemma}

\begin{proof}
Given a nonempty pattern $v$, let us substitute $a$ for each occurrence of
$*$ in $v$ and $b^{-1}$ for each occurrence of $*^{-1}$.  We get a word
$\xi_1\in Q_\pm^*$ that follows the pattern $v$.  Now let us modify $\xi_1$
by changing its first letter.  If this letter is $a$, we change it to $c$.
If the first letter of $\xi_1$ is $b^{-1}$, we change it to $c^{-1}$.  This
yields another word $\xi_2\in Q_\pm^*$ that follows the pattern $v$.  By
construction, $\xi_1$ and $\xi_2$ are freely irreducible.  Since $\chi(a)=
\chi(b^{-1})=-1$ and $\chi(c)=\chi(c^{-1})=1$, it follows that
$\chi(\xi_2)=-\chi(\xi_1)$.
\end{proof}

The following proposition is the main step in the proof of Theorem
\ref{main1}.

\begin{proposition}\label{ind3}
Suppose $\xi\in Q_\pm^*$ is a freely irreducible word.  Then the orbit of
$\xi$ under the action of the group $G(E)$ on $Q_\pm^*$ consists of all
freely irreducible words following the same pattern as $\xi$.
\end{proposition}

To prove Proposition \ref{ind3}, we need several lemmas.

\begin{lemma}\label{ind4}
Two words $\xi_1,\xi_2\in Q_\pm^*$ are in the same orbit of the $G(E)$
action if and only if the reversed words $\rev{\xi_1}$ and $\rev{\xi_2}$
are in the same orbit of this action.
\end{lemma}

\begin{proof}
Using Proposition \ref{auto3}, we verify that the reverse automaton $R$ of
$E$ is well defined.  The Moore diagram $\Gamma_R$ of $R$ is obtained from
the Moore diagram $\Gamma_E$ of $E$ by reversing all edges.
It is easy to observe that $\Gamma_R$ can also be obtained by renaming
vertices $\alpha$ and $\beta$ of the graph $\Gamma_E$ to $\beta$ and
$\alpha$, respectively (see Figure \ref{fig4}).  Hence $R_\alpha=E_\beta$,
$R_\beta=E_\alpha$, $R_\gamma=E_\gamma$.  In particular, $S(R)=S(E)$.  By
Lemma \ref{free2}, $E_\alpha$, $E_\beta$, and $E_\gamma$ are involutions;
therefore $S(E)=G(E)$.  Now it follows from Proposition \ref{auto3} that,
given $\xi_1,\xi_2\in Q_\pm^*$, we have $\xi_2=g(\xi_1)$ for some $g\in
G(E)$ if and only if $\rev{\xi_2}=h(\rev{\xi_1})$ for some $h\in G(E)$.
\end{proof}

Let $V^0$ denote the set of patterns without double letters.  That is, a
pattern $v\in V^0$ does not contain subwords $**$ and $*^{-1}*^{-1}$.  We
consider $4$ subsets of $V^0$.  By $V_{++}^0$ denote the set of $v\in V^0$
such that $*$ is the first and the last letter of $v$.  By $V_{+-}^0$
denote the set of $v\in V^0$ such that either $v=\varnothing$ or the first
letter of $v$ is $*$ while the last letter is $*^{-1}$.  By $V_{-+}^0$
denote the set of $v\in V^0$ such that either $v=\varnothing$ or the first
letter of $v$ is $*^{-1}$ while the last letter is $*$.  By $V_{--}^0$
denote the set of $v\in V^0$ such that $*^{-1}$ is the first and the last
letter of $v$.  Clearly, $V^0=V_{++}^0\cup V_{+-}^0\cup V_{-+}^0\cup
V_{--}^0$.

Now we define $W_{++}^0$ (resp. $W_{+-}^0$, $W_{-+}^0$, $W_{--}^0$) as the
set of all words over the alphabet $\{a,b,a^{-1},b^{-1}\}$ that follow
patterns from the set $V_{++}^0$ (resp. $V_{+-}^0$, $V_{-+}^0$,
$V_{--}^0$).  In particular, $W_{++}^0$ contains the words $a,ab^{-1}a,
ab^{-1}ab^{-1}a,\dots$, $W_{+-}^0$ contains $\varnothing,ab^{-1},
ab^{-1}ab^{-1},\dots$, $W_{-+}^0$ contains $\varnothing,b^{-1}a,
b^{-1}ab^{-1}a,\dots$, and $W_{--}^0$ contains $b^{-1},b^{-1}ab^{-1},
b^{-1}ab^{-1}ab^{-1},\dots$.

Consider an endomorphism $r$ of the free monoid $Q_\pm^*$ defined by
$r(a)=a$, $r(a^{-1})=a^{-1}$, $r(b)=b$, $r(b^{-1})=b^{-1}$, $r(c)=
r(c^{-1})=\varnothing$.  For any $\xi\in Q_\pm^*$ the word $r(\xi)$ is
obtained by deleting all letters $c$ and $c^{-1}$ in $\xi$.  The
restriction of $r$ to the set $W_{++}^0\cup W_{+-}^0\cup W_{-+}^0\cup
W_{--}^0$ is the identity map.  Let $W_{++}=r^{-1}(W_{++}^0)$,
$W_{+-}=r^{-1}(W_{+-}^0)$, $W_{-+}=r^{-1}(W_{-+}^0)$, and
$W_{--}=r^{-1}( W_{--}^0)$.

\begin{lemma}\label{indextra}
(i) $E_\alpha(\xi)=\xi$ for all $\xi\in W_{++}\cup W_{+-}$ while
$E_\beta(\xi)=\xi$ for all $\xi\in W_{-+}\cup W_{--}$.

(ii) If $\xi\in W_{++}$ then $E_\alpha(\xi a)=\xi b$ and $E_\alpha(\xi b)=
\xi a$.  If $\xi\in W_{+-}$ then $E_\alpha(\xi a^{-1})=\xi b^{-1}$ and
$E_\alpha(\xi b^{-1})=\xi a^{-1}$.  If $\xi\in W_{-+}$ then $E_\beta(\xi
a)=\xi b$ and $E_\beta(\xi b)=\xi a$.  If $\xi\in W_{--}$ then
$E_\beta(\xi a^{-1})=\xi b^{-1}$ and $E_\beta(\xi b^{-1})=\xi a^{-1}$.
\end{lemma}

\begin{proof}
Let $\xi\in Q_\pm^*$.  If $x\in\{c,c^{-1}\}$ then $E_\alpha(x\xi)=
xE_\alpha(\xi)$ and $E_\beta(x\xi)=xE_\beta(\xi)$.  Besides, $\xi$ belongs
to one of the sets $W_{++}$, $W_{+-}$, $W_{-+}$, $W_{--}$ if and only if
$x\xi$ belongs to the same set.  Further, if $x\in\{a,b\}$ then
$E_\alpha(x\xi)=xE_\beta(\xi)$.  In this case, $x\xi\in W_{++}\cup W_{+-}$
if and only if $\xi\in W_{-+}\cup W_{--}$.  Finally, if $x\in
\{a^{-1},b^{-1}\}$ then $E_\beta(x\xi)=xE_\alpha(\xi)$.  In this case,
$x\xi\in W_{-+}\cup W_{--}$ if and only if $\xi\in W_{++}\cup W_{+-}$.  The
statement (i) of the lemma follows from the above by induction on the
length of $\xi$.

Any word $\xi\in W_{++}$ contains an odd number of letters $a,b,a^{-1},
b^{-1}$.  Therefore $E_\alpha(\xi\eta)=E_\alpha(\xi)E_\beta(\eta)=\xi
E_\beta(\eta)$ for all $\eta\in Q_\pm^*$.  In particular, $E_\alpha(\xi a)=
\xi b$ and $E_\alpha(\xi b)=\xi a$.  The rest of the statement (ii) is
obtained in a similar way.  We omit the details.
\end{proof}

Given a freely irreducible word $\xi\in Q_\pm^*$, let $Z(\xi)$ denote the
set of freely irreducible words in $Q_\pm^*$ that follow the same pattern
as $\xi$ and match $\xi$ completely or except for the last letter.
Obviously, $\xi\in Z(\xi)$, and $\eta\in Z(\xi)$ if and only if $\xi\in
Z(\eta)$.  If $\xi\ne\varnothing$ then $Z(\xi)$ consists of $2$ or $3$
words.  Namely, there are exactly $3$ words in $Q_\pm^*$ that follow the
same pattern as $\xi$ and match $\xi$ completely or except for the last
letter.  However if the last two letters in the pattern of $\xi$ are
distinct then one of these $3$ words is freely reducible.

\begin{lemma}\label{ind5}
Suppose $v$ is a pattern of length at least $2$ such that the last two
letters of $v$ are distinct.  Then there exist two freely irreducible words
$\xi_a,\xi_b\in Q_\pm^*$ such that $\xi_a$ and $\xi_b$ follow the pattern
$v$, $Z(\xi_a)=\{\xi_a,\xi_b\}$, and $\xi_b=g(\xi_a)$, where
$g\in\{E_\alpha,E_\beta\}$.
\end{lemma}

\begin{proof}
Let $v_0$ denote the pattern obtained by deleting the last letter of $v$.
We replace each letter $*$ in the word $v_0$ by $a$ if the next letter is
$*^{-1}$ and by $c$ otherwise.  Each letter $*^{-1}$ is replaced by
$b^{-1}$ if the next letter is $*$ and by $c^{-1}$ otherwise.  We get a
word $\eta\in Q_\pm^*$ that follows the pattern $v_0$.  For example, if
$v_0=****^{-1}*^{-1}*^{-1}****^{-1}*^{-1}*^{-1}$ then
$\eta=ccac^{-1}c^{-1}b^{-1}ccac^{-1}c^{-1}c^{-1}$.  Now let $\xi_a=\eta a$,
$\xi_b=\eta b$ if $v=v_0*$ and let $\xi_a=\eta a^{-1}$, $\xi_b=\eta b^{-1}$
if $v=v_0*^{-1}$.  The words $\xi_a$ and $\xi_b$ follow the pattern $v$.
By construction, $\eta$ is freely irreducible and its last letter is $c$ or
$c^{-1}$.  Since the last letters of patterns $v$ and $v_0$ are different,
it follows that $\xi_a$, $\xi_b$ are freely irreducible and $Z(\xi_a)=
\{\xi_a,\xi_b\}$.

It is easy to observe that $\eta$ belongs to one of the sets $W_{++}$,
$W_{+-}$, $W_{-+}$, $W_{--}$.  Lemma \ref{indextra} implies that $g(\xi_a)=
\xi_b$ and $g(\xi_b)=\xi_a$, where $g\in\{E_\alpha,E_\beta\}$.  Namely,
$g=E_\alpha$ if the first letter of $v$ is $*$ and $g=E_\beta$ otherwise.
\end{proof}

\begin{lemma}\label{ind6}
Suppose $v$ is a pattern of length at least $2$ such that the first two
letters of $v$ coincide as well as the last two letters of $v$.  Then there
exist $9$ freely irreducible words $\xi_{q_1q_2}\in Q_\pm^*$, $q_1,q_2\in
\{a,b,c\}$ such that each $\xi_{q_1q_2}$ follow the pattern $v$ and for any
$q\in\{a,b,c\}$ we have $Z(\xi_{qc})=\{\xi_{qa},\xi_{qb},\xi_{qc}\}$,
$Z(\rev{\xi_{cq}})=\{\rev{\xi_{aq}},\rev{\xi_{bq}},\rev{\xi_{cq}}\}$, and
$\xi_{qb}=g_q(\xi_{qa})$, where $g_q\in\{E_\alpha,E_\beta\}$.
\end{lemma}

\begin{proof}
Let $v_0$ denote the pattern obtained by deleting the first and the last
letters of $v$.  We replace each letter $*$ in the word $v_0$ by $a$ if the
previous letter is $*^{-1}$ and by $c$ otherwise.  Each letter $*^{-1}$ is
replaced by $b^{-1}$ if the previous letter is $*$ and by $c^{-1}$
otherwise.  We get a word $\eta\in Q_\pm^*$ that follows the pattern $v_0$.
For example, if $v_0=****^{-1}*^{-1}*^{-1}****^{-1}*^{-1}*^{-1}$ then
$\eta=cccb^{-1}c^{-1}c^{-1}accb^{-1}c^{-1}c^{-1}$.  Now for any
$q\in\{a,b,c\}$ let $\xi_q=q\eta$ if the first letter of $v$ is $*$ and let
$\xi_q=q^{-1}\eta$ otherwise.  Further, for any $q_1,q_2\in\{a,b,c\}$ let
$\xi_{q_1q_2}=\xi_{q_1}q_2$ if the last letter of $v$ is $*$ and let
$\xi_{q_1q_2}=\xi_{q_1}q_2^{-1}$ otherwise.  Each $\xi_{q_1q_2}$ follows
the pattern $v$.  By construction, $\eta$ is freely irreducible.  Since the
first two letters of $v$ coincide and so do the last two letters, it
follows that each $\xi_{q_1q_2}$ is freely irreducible.  Then
$Z(\xi_{qc})=\{\xi_{qa},\xi_{qb},\xi_{qc}\}$ and $Z(\rev{\xi_{cq}})=
\{\rev{\xi_{aq}},\rev{\xi_{bq}},\rev{\xi_{cq}}\}$ for any $q\in\{a,b,c\}$.

It is easy to observe that the words $\xi_a$ and $\xi_b$ belong to one of
the sets $W_{++}\cup W_{+-}$ and $W_{-+}\cup W_{--}$ while $\eta$ and
$\xi_c$ belong to the other.  Lemma \ref{indextra} implies that
$g(\xi_{aa})=\xi_{ab}$, $g(\xi_{ab})=\xi_{aa}$, $g(\xi_{ba})=\xi_{bb}$, and
$g(\xi_{bb})=\xi_{ba}$, where $g=E_\alpha$ if the first letter of $v$ is
$*$ and $g=E_\beta$ otherwise.  Moreover, $h(\xi_{ca})=\xi_{cb}$ and
$h(\xi_{cb})=\xi_{ca}$, where $h$ is the element of $\{E_\alpha,E_\beta\}$
different from $g$.
\end{proof}

\begin{proofof}{Proposition \ref{ind3}}
First we shall show that the $G(E)$ action on $Q_\pm^*$ preserves patterns
and free irreducibility of words.  Suppose $q_1|q_2$ is the label assigned
to an edge of the Moore diagram of the automaton $E$.  Then either
$q_1,q_2\in\{a,b,c\}$ or $q_1,q_2\in\{a^{-1},b^{-1},c^{-1}\}$ (see Figure
\ref{fig4}).  It follows that transformations $E_\alpha$, $E_\beta$, and
$E_\gamma$ preserve patterns of words.  So does any $g\in G(E)$.  Further,
it is easy to verify that sets $P_1=\{aa^{-1},bb^{-1},cc^{-1}\}$ and
$P_2=\{a^{-1}a,b^{-1}b,c^{-1}c\}$ are invariant under transformations
$E_\alpha$, $E_\beta$, and $E_\gamma$.  Therefore $P_1$ and $P_2$ are
invariant under the $G(E)$ action on $Q_\pm^*$.  Any freely reducible word
$\xi\in Q_\pm^*$ is represented as $\xi_1\xi_0\xi_2$, where $\xi_0\in
P_1\cup P_2$ and $\xi_1,\xi_2\in Q_\pm^*$.  Proposition \ref{auto4} implies
that for any $g\in S(E)=G(E)$ we have $g(\xi)=g(\xi_1)g_0(\xi_0)
g_1(\xi_2)$, where $g_0,g_1\in G(E)$.  By the above $g(\xi)$ is freely
reducible.  Thus the $G(E)$ action preserves free reducibility of words.
Since $G(E)$ is a group, its action on $Q_\pm^*$ preserves free
irreducibility as well.

Now we are going to prove that for any freely irreducible words
$\xi_1,\xi_2\in Q_\pm^*$ following the same pattern $v$ there exists
$g\in G(E)$ such that $\xi_2=g(\xi_1)$.  The claim is proved by induction
on the length of the pattern $v$.  The empty pattern is followed only by
the empty word.  As for one-letter patterns, it is enough to notice that
$\pi_{(abc)}(a)=b$, $\pi_{(abc)}(b)=c$, $\pi_{(abc)}(c)=a$,
$\pi_{(abc)}(a^{-1})=b^{-1}$, $\pi_{(abc)}(b^{-1})=c^{-1}$,
$\pi_{(abc)}(c^{-1})=a^{-1}$, and $\pi_{(abc)}\in G(E)$ by Lemma
\ref{free2}.  Now let $n\ge2$ and assume that the claim holds for all
patterns of length less than $n$.  Take any pattern $v$ of length $n$.
First consider the case when the last two letters of $v$ are distinct.  By
Lemma \ref{ind5}, there exists a freely irreducible word $\xi\in Q_\pm^*$
such that $\xi$ follows the pattern $v$ and the set $Z(\xi)$ is contained
in an orbit of the $G(E)$ action.  Secondly, consider the case when the
first two letters of $v$ are distinct.  Here, the last two letters of the
reversed pattern $\rev{v}$ are distinct.  Then there is a freely
irreducible word $\xi_-$ such that $\xi_-$ follows the pattern $\rev{v}$
and $Z(\xi_-)$ is contained in an orbit of the $G(E)$ action.  Clearly,
$\xi_-=\rev{\xi_+}$, where $\xi_+$ is a freely irreducible word following
the pattern $v$.

Now consider the case when the last two letters of $v$ coincide as well as
the first two letters.  Let $\xi_{q_1q_2}$, $q_1,q_2\in\{a,b,c\}$ be the
words provided by Lemma \ref{ind6}.  By $\xi_a$, $\xi_b$, $\xi_c$ denote
the words obtained by deleting the last letter of $\xi_{aa}$, $\xi_{ba}$,
$\xi_{ca}$, respectively.  Obviously, $\xi_a$, $\xi_b$, and $\xi_c$ are
freely irreducible and follow the same pattern of length $n-1$.  By the
inductive assumption there exist $h_a,h_b\in G(E)$ such that
$h_a(\xi_c)=\xi_a$, $h_b(\xi_c)=\xi_b$.  Since the $G(E)$ action preserves
patterns, it follows that $h_a(\xi_{cc})\in Z(\xi_{ac})$, $h_b(\xi_{cc})\in
Z(\xi_{bc})$.  Recall that there are involutions $g_a,g_b,g_c\in G(E)$ such
that $g_q(\xi_{qa})=\xi_{qb}$, $q\in\{a,b,c\}$.  If $h_a(\xi_{cc})=
\xi_{aa}$ or $h_a(\xi_{cc})=\xi_{ab}$, then the word
$h_a^{-1}g_ah_a(\xi_{cc})$ matches $\xi_{cc}$ except for the last letter.
Hence $h_a^{-1}g_ah_a(\xi_{cc})\in\{\xi_{ca},\xi_{cb}\}$.  Since
$g_c(\xi_{ca})=\xi_{cb}$, it follows that $Z(\xi_{cc})=
\{\xi_{ca},\xi_{cb},\xi_{cc}\}$ is contained in an orbit of the $G(E)$
action.  If $h_b(\xi_{cc})=\xi_{ba}$ or $h_b(\xi_{cc})=\xi_{bb}$, we reach
the same conclusion.  On the other hand, if $h_a(\xi_{cc})=\xi_{ac}$ and
$h_b(\xi_{cc})=\xi_{bc}$ then the set $\{\xi_{ac},\xi_{bc},\xi_{cc}\}$ is
contained in an orbit.  By Lemma \ref{ind4}, $Z(\rev{\xi_{cc}})=
\{\rev{\xi_{ac}},\rev{\xi_{bc}},\rev{\xi_{cc}}\}$ is also contained in an
orbit of the $G(E)$ action.

In any of the considered cases, there exists a freely irreducible word
$\xi$ such that $\xi$ follows the pattern $v$ and one of the sets $Z(\xi)$,
$Z(\rev{\xi})$ is contained in an orbit of the $G(E)$ action.  Take any
freely irreducible words $\xi_1$, $\xi_2$ following the pattern $v$.
First suppose that $Z(\xi)$ is contained in an orbit.  Let $\eta,\eta_1,
\eta_2$ be words obtained by deleting the last letter of $\xi,\xi_1,\xi_2$,
respectively.  Then $\eta,\eta_1,\eta_2$ are freely irreducible and follow
the same pattern obtained by deleting the last letter of $v$.  By the
inductive assumption there are $g_1,g_2\in G(E)$ such that
$\eta=g_1(\eta_1)=g_2(\eta_2)$.  Since the $G(E)$ action preserves patterns
and free irreducibility, it follows that $g_1(\xi_1),g_2(\xi_2)\in Z(\xi)$.
As $Z(\xi)$ is contained in an orbit, there exists $g\in G(E)$ such that
$g(g_1(\xi_1))=g_2(\xi_2)$.  Then $\xi_2=g_2^{-1}gg_1(\xi_1)$.
Now suppose that the set $Z(\rev{\xi})$ is contained in an orbit.  Let
$\zeta,\zeta_1,\zeta_2$ be words obtained by deleting the last letter of
$\rev{\xi},\rev{\xi_1},\rev{\xi_2}$, respectively.  Clearly,
$\zeta,\zeta_1,\zeta_2$ are freely irreducible and follow the same pattern
obtained by deleting the last letter of $\rev{v}$.  By the inductive
assumption there exist $h_1,h_2\in G(E)$ such that $\zeta=h_1(\zeta_1)=
h_2(\zeta_2)$.  Since the $G(E)$ action preserves patterns and free
irreducibility, we conclude that $h_1(\rev{\xi_1}),h_2(\rev{\xi_2})\in
Z(\rev{\xi})$.  It follows that $\rev{\xi_1}$ and $\rev{\xi_2}$ are in the
same orbit of the $G(E)$ action.  By Lemma \ref{ind4}, $\xi_1$ and $\xi_2$
are also in the same orbit.  The induction step is complete.
\end{proofof}

\begin{proofof}{Theorem \ref{main1}}
The group $G(A)$ is the free non-Abelian group on generators $A_a$, $A_b$,
$A_c$ if and only if $A_{q_1}^{m_1}A_{q_2}^{m_2}\dots A_{q_n}^{m_n}\ne1$
for any pair of sequences $q_1,\dots,q_n$ and $m_1,\dots,m_n$ such that
$n>0$, $q_i\in\{a,b,c\}$ and $m_i\in\bZ\setminus\{0\}$ for $1\le i\le n$,
and $q_i\ne q_{i+1}$ for $1\le i\le n-1$.  Since $B_a=A_a$, $B_b=A_b$,
$B_c=A_c$, $B_{a^{-1}}=A_a^{-1}$, $B_{b^{-1}}=A_b^{-1}$, $B_{c^{-1}}=
A_c^{-1}$, an equivalent condition is that $B_\xi\ne1$ for any nonempty
freely irreducible word $\xi\in Q_\pm^*$.

Suppose $B_\xi=1$ for some freely irreducible word $\xi\in Q_\pm^*$.  By
Corollary \ref{auto5}, $B_{g(\xi)}=1$ for all $g\in S(D)$.  Then
Proposition \ref{auto1} imply that $B_{g(\xi)}=1$ for all $g\in G(D)$.  By
Proposition \ref{free4}, $G(D)=G(E)$.  Now it follows from Proposition
\ref{ind3} that $B_\eta=1$ for any freely irreducible word $\eta\in
Q_\pm^*$ following the same pattern as $\xi$.  In particular, $B_\eta$ acts
trivially on the first level of the binary tree $\{0,1\}^*$.  Finally,
Lemmas \ref{ind1} and \ref{ind2} imply that $\xi$ follows the empty
pattern.  Then $\xi$ itself is the empty word.
\end{proofof}

\medskip

{\sc
\begin{raggedright}
Department of Mathematics\\
Texas A\&M University\\
College Station, TX 77843--3368
\end{raggedright}
}

\end{document}